\documentclass[12pt]{article}
\usepackage[latin1]{inputenc}
\usepackage{amssymb}
\newtheorem{theorem}{Theorem}

\newtheorem{proposition}[theorem]{Proposition}
\newtheorem{sublemma}[theorem]{Sublemma}
\newtheorem{lemma}[theorem]{Lemma}
\newtheorem{corollary}[theorem]{Corollary}
\newtheorem{remark}[theorem]{Remark}

\newcommand{\R}{\mathbb{R}}
\newcommand{\Q}{\mathbb{Q}}

\newcommand{\Sf}{\mathbb{S}}

\newcommand{\Les}{\mathbb{L}}

\newcommand{\C}{\mathbb{C}}
\newcommand{\Hy}{\mathbb{H}}

\newcommand{\spa}{\mbox{span}}
\newcommand{\hess}{\mbox{Hess\,}}

\newcommand{\po}{{\hspace*{-1ex}}{\bf .  }}

\def\<{\langle}

\def\>{\rangle}

\def\bea{\begin{eqnarray*} }
\def\eea{\end{eqnarray*} }
\def\be{\begin{equation} }
\def\ee{\end{equation} }

\def\proof{\noindent{\it Proof: }}
\def\qed{\ifhmode\unskip\nobreak\fi\ifmmode\ifinner
\else\hskip5 pt \fi\fi\hbox{\hskip5 pt \vrule width4 pt
height6 pt  depth1.5 pt \hskip 1pt }}
\begin{document}
\title{\bigskip
\bigskip
Submanifolds of codimension two attaining  equality\\ in an extrinsic inequality}
\author{Marcos Dajczer $\&$ Ruy Tojeiro}
\date{}
\maketitle

\begin{abstract}
We provide a parametric construction in terms of minimal surfaces of the Euclidean submanifolds of  codimension two and arbitrary  dimension that attain equality in an inequality due to De Smet, Dillen, Verstraelen and Vrancken. The latter  involves the scalar curvature, the norm of the normal curvature tensor and the length of the mean curvature vector.  
\end{abstract}

Let $f\colon\,M^n\to \Q_c^{n+p}$ be an isometric immersion of an $n$-dimensional  Riemannian \mbox{manifold} into a space form of dimension $n+p$ and constant sectional curvature $c$. Let $s$ denote the  normalized scalar curvature of $M^n$ and let $s_N$ be given by $$n(n-1)s_N=\|R^\perp\|,$$ where $R^\perp$ is the normal curvature tensor of $f$. Explicitly, 
$$
s=\frac{2}{n(n-1)}\sum_{1\leq i<j\leq n}\<R(e_i,e_j)e_j,e_i\>
$$
and 
$$
s_N=\frac{2}{n(n-1)}(\!\!\sum_{1\leq i<j\leq n
\atop 1\leq r<s\leq p} \!\!\<R^\perp(e_i,e_j)\xi_r,\xi_s\>^2)^{1/2},
$$
 where $R$  is the curvature tensor of $M^n$ and $\{e_1,\ldots, e_n\}$ (resp.,  $\{\xi_1,\ldots,\xi_p\}$) is an orthonormal basis of the tangent (resp., normal) space.   
%Observe that $n(n-1)s_N$ is just the norm of the tensor $R^\perp$.

 The pointwise inequality
$$
\hspace{33ex} 
s\leq c+\|H\|^2-s_N
%s+s_N\leq \|H\|^2+c
\hspace{30ex}(*)
$$
relates the intrinsic scalar curvature $s$ to the extrinsic data on the right-hand-side. Here $H$ denotes  the mean curvature vector  of $f$. It was proved for codimension $p=2$ by De Smet, Dillen, Verstraelen and Vrancken in \cite{ddvv}. Also, the pointwise structure of the shape operators of submanifolds attaining equality  was determined. It was shown that  equality holds at $x\in M^n$ if and only if there exist orthonormal bases $\{e_1,\ldots, e_n\}$ and $\{\eta, \zeta\}$ of the tangent and normal spaces at $x$, respectively,  such that the shape operators $A_\eta$ and $A_\zeta$ have the  form

\be\label{eq:shapeop}
A_\eta=\left[ 
\begin{array}{ccc}
\lambda & \;\mu &\;0 \;\;\cdots\;\; 0\\ 
\mu & \lambda &\;0 \;\;\cdots\;\; 0\\ 
0 & 0 &\;\lambda \;\;\cdots\;\; 0\\ 
\vdots & \vdots & \vdots \;\;\cdots\;\; \vdots\\ 
0 & 0 &\;0 \;\;\cdots\;\; \lambda\\
\end{array}
\right];\;\;\;\;\;\;\;\; A_\zeta=\left[ 
\begin{array}{ccc}
\mu & \;0 &\;0 \;\;\cdots\;\; 0\\ 
0 & -\mu &\;0 \;\;\cdots\;\; 0\\ 
0 & 0 &\;0 \;\;\cdots\;\; 0\\ 
\vdots & \vdots & \vdots \;\;\cdots\;\; \vdots\\ 
0 & 0 &\;0 \;\;\cdots\;\; 0\\
\end{array}
\right].\ee

    The  inequality $(*)$ is also known to hold for surfaces \cite{gr} and for submanifolds with flat normal bundle  \cite{ch} of any codimension 
%(in particular for hypersurfaces) 
as well as for various special classes of submanifolds (see \cite{dfv} and the references therein). Moreover, it was conjectured in \cite{ddvv} to hold for any submanifold of a space form. 
Recently, the conjecture was proved in \cite{cl} and \cite{cl2}, respectively,    for  three dimensional submanifolds with arbitrary codimension and for   submanifolds of codimension three and any codimension of space forms.

For any isometric immersion  $f\colon\,M^n\to \Q_c^{n+p}$, it was shown in \cite{dfv} that $(*)$
holds at a point $x\in M^n$ if and only if the inequality 
$$
\sum_{\alpha, \beta=1}^p\|[B_\alpha, B_\beta]\|^2\leq (\sum_{\alpha=1}^p\|B_\alpha\|^2)^2
$$
is satisfied for the traceless parts $B_1,\ldots, B_p$ of the shape operators of $f$ with respect to any orthonormal normal frame at $x$.

An important consequence of this reformulation of the  inequality is that it readily implies that the class of submanifolds $f\colon\,M^n\to \Q_c^{n+p}$ for which  equality holds 
is invariant under conformal transformations of the ambient space. In fact, under such a transformation the traceless parts of the shape operators only change by multiplication by a common smooth function on $M^n$. 

By the above, the class of isometric immersions $f\colon\,M^n\to \R^{n+2}$
attaining  equality everywhere  in  $(*)$ contains any composition
of an inversion in  $\R^{n+2}$ with  a minimal isometric immersion
$f\colon\,M^n\to \R^{n+2}$ whose shape operators are as in (\ref{eq:shapeop}) with $\lambda=0$.  Notice that such minimal submanifolds belong to the class of  austere submanifolds of rank two, first studied in arbitrary codimension by Bryant \cite{br} for $n=3$ and then by Dajczer-Florit \cite{df2} for any dimension $n$.

   In this paper, we provide an explicit local construction of all Euclidean submanifolds  $f\colon\,M^n\to \R^{n+2}$ attaining equality   everywhere  in the inequality $(*)$. For $n=2$, such submanifolds are precisely the surfaces in $\R^4$  whose ellipses of curvature are circles at any point, and this was considered in our previous paper \cite{dt}. It turns out that several steps of the proof of the main result of that paper can be adapted to the present general case.

 Our construction starts with a simply connected minimal surface
$g\colon\,M^2\to \R^{n+2}$, oriented by a global conformal
diffeomorphism  onto either the complex plane or the unit disk.
Then we consider its conjugate  minimal surface $h\colon\,M^2\to
\R^{n+2}$, each of whose components with respect to this global
parameter is the harmonic conjugate of the corresponding component
of $g$. Equivalently, $h_*=g_*\circ J$,
where $J$ is the complex structure on $M^2$ compatible with its
orientation. Now we decompose the  position vector of $h$ in its tangent and normal components with respect to $g$, i.e.,
$$
h=g_*h^T + h^N.
$$
Finally, on the complement of the subset of isolated points of $M$ where $h^N$ vanishes, let $\Lambda_1$ be the unit bundle of the vector subbundle $\Lambda$  of the normal bundle  of $g$ that is orthogonal to $h^N$. We can now state our main result.

\begin{theorem}\po\label{thm:main} Assume that $n\geq 3$ and  define a map $\phi\colon\,\Lambda_1\to\R^{n+2}$ by
\be\label{eq:phi}
\phi(y,w)=g(y)+g_*Jh^T(y)+\|h^N(y)\|w.
\ee
 Then, at regular points, $\phi$ parameterizes an $n$-dimensional  submanifold in $\R^{n+2}$ attaining  equality in the inequality $(*)$. 

Conversely, any submanifold $f\colon\,M^n\to\R^{n+2}$, $n\geq 3$,  free of umbilical and minimal points and attaining  equality in  the inequality $(*)$ can be parameterized in this way.
\end{theorem}

By combining the preceding result with the generalized Weierstrass
parameterization of Euclidean minimal surfaces $g\colon\,M^2\to \R^{n+2}$ (cf.\ \cite{ho},\,\cite{la}) we have a parametric
representation of the submanifolds $f\colon\,M^n\to\R^{n+2}$  attaining  equality in  $(*)$.
\medskip

We also characterize in terms of our construction the submanifolds that are images  of  austere submanifolds of rank two attaining equality in the inequality $(*)$ in either $\mathbb{R}^{n+2}$, $\mathbb{S}_d^{n+2}$ or $\mathbb{H}_d^{n+2}$ by an inversion in the first case or a stereographic projection in the other two cases.  Here, \mbox{$\mathbb{S}_d^{n+2}=\mathbb{S}^{n+2}(de_{n+3};d)$} is
the sphere in $\mathbb{R}^{n+3}$ with radius $d$ centered at $de_{n+3}$ and
$\mathbb{H}_d^{n+2}=\mathbb{H}^{n+2}(-de_{n+3};d)$ is the hyperbolic space
$$
\mathbb{H}_d^{n+2}=\{X\in \mathbb{L}^{n+3}\,:\,\langle X+de_{n+3},X+de_{n+3}\rangle=-d^2\}.
$$
Moreover, we regard $\R^{n+2}$ as the
hyperplane through the origin and normal to the unit vector $e_{n+3}$ in either $\mathbb{R}^{n+3}$ or
Lorentzian space $\mathbb{L}^{n+3}$, and by the stereographic projection of $\Hy_d^{n+2}$ onto the open ball $B(0;2d)\subset \R^{n+2}$ we mean the map that assigns to each $P\in \Hy_d^{n+2}$ 
the point of $\R^{n+2}$ where the line through the points  $-2de_{n+3}$ and $P$ intersects $\mathbb{R}^{n+2}$.  
Let
$$
G=g+ih\colon\,L^2\to\C^{n+2}\approx\R^{n+2}+i\R^{n+2}
$$
be the  holomorphic representative of the minimal surface $g\colon\,L^2\to \R^{n+2}$.
For any real number $k$, we denote by ${\cal H}_k^{n+1}\subset\C^{n+2}$  the quadric 
$$
{\cal H}^{n+1}_k=\{Z\in\C^{n+2}: \<\!\< Z,Z\>\!\>=k\},
$$
where $\<\!\< \,\,,\,\,\>\!\>\colon\,\C^{n+2}\times \C^{n+2}\to \C$  denotes the  linear inner
product on
$\C^{n+2}$.

\begin{theorem}\po\label{prop:minimas} 
The map $\phi$ in (\ref{eq:phi}) parameterizes the composition of an austere submanifold of rank two in $\R^{n+2}$, $\Sf_d^{n+2}$ or $\Hy_d^{n+2}$ attaining equality in the inequality $(*)$ with an inversion in $\R^{n+2}$ with respect to a hypersphere centered at the origin or a stereographic projection of $\Sf_d^{n+2}$ or $\Hy_d^{n+2}$ onto $\R^{n+2}$ and $B(0;2d)\subset \R^{n+2}$, respectively, if and only if 
%the holomorphic representative 
$G$ 
%of $g$ 
takes values  in ${\cal H}_k^{n+1}$, with $k=0,\,4d^2$ or $-4d^2$,  respectively.
\end{theorem}

As a consequence of Theorems \ref{thm:main} and \ref{prop:minimas}, we obtain the following parameterization of all  austere $n$-dimensional submanifolds of 
rank two in $\R^{n+2}$ attaining equality in $(*)$.

\begin{corollary}\po\label{cor:austere} Any austere $n$-dimensional submanifold of rank two  in $\R^{n+2}$, $\Sf_d^{n+2}$ or $\Hy_d^{n+2}$ attaining equality in the inequality $(*)$ can be parameterized~by
$$
\psi={\cal I}\circ \phi,
$$
where $\phi$ is given by (\ref{eq:phi}) in terms of a minimal surface $g\colon\,L^2\to \R^{n+2}$ whose holomorphic representative $g+ih$ takes values in a quadric ${\cal H}_k^{n+1}$ of $\C^{n+2}$ with $k=0,\,4d^2$ or $-4d^2$,  respectively, and ${\cal I}$ is an inversion with respect to a hypersphere centered at the origin or the inverse of a stereographic projection of $\Sf_d^{n+2}$ or $\Hy_d^{n+2}$ onto $\R^{n+2}$ or $B(0;2d)\subset \R^{n+2}$, respectively.
\end{corollary}

Our next result extends  Theorem $2$ in \cite{dt}. Let $
G=g+ih\colon\,L^2\to\C^{n+2}$ be the  holomorphic
representative 
of  the minimal surface
$g\colon\,L^2\to \R^{n+2}$ associated to a submanifold $\phi\colon\,M^n\to\R^{n+2}$ attaining equality in the inequality $(*)$. Let   \mbox{$\tilde
G=\tilde{g}+i\tilde{h}$} be the holomorphic representative of  the minimal surface
$\tilde{g}\colon\,L^2\to \R^{n+2}$  associated to its composition
$\tilde{\phi}={\cal I}\circ \phi$   with an
inversion ${\cal I}$ in $\R^{n+2}$ with respect to a sphere of
radius $d$ taken, for simplicity, centered at the origin.  Then $G$ and $\tilde{G}$ are related as follows.

\begin{theorem}\po\label{thm:pairs}
If $\phi$ is not the composition of an austere submanifold attaining equality in the inequality $(*)$ with an inversion, then $\tilde G= \overline{T_d \circ G}$, where $T_d=d^{\,2}T$ and $T\colon\,\C^{n+2}\to \C^{n+2}$ is the
holomorphic map
$$
T(Z)=\frac{Z}{\<\!\< Z,Z\>\!\>}.
$$
\end{theorem}

  As pointed out in \cite{dt},  the holomorphic map $T_d\colon\,\C^{m}\to \C^{m}$ for any $m$ and $T$ defined above,  can be regarded as the {\em inversion\/} in  $\C^{m}$ with respect to the quadric ${\cal H}^{m-1}_{d^2}$.
\medskip  

  The following  byproduct of Theorem \ref{thm:pairs} yields,  in particular, a transformation for minimal surfaces in $\R^{n+2}$.
\begin{corollary}\po\label{cor:pairs}
 The holomorphic inversion map
$T$ preserves the class of holomorphic curves
$G=g+ih\colon\,L^2\to \C^{n+2}$ whose real and imaginary parts $g$
and $h$ define  conjugate minimal immersions into $\R^{n+2}$.
\end{corollary}

   After the paper was submitted for publication, two independent proofs of  the inequality $(*)$ for submanifolds of arbitrary codimension have appeared (see   \cite{gt} and \cite{cl3}). Moreover, in \cite{gt} also the pointwise structure of the second fundamental forms of the submanifolds that attain equality was determined. In particular, it was shown that the  first normal spaces of such submanifolds, i.e., the subspaces of the normal spaces that are spanned by the image of the second fundamental form, have dimension  either two or three. If the first case holds everywhere and the submanifold has dimension at least four and is not minimal, then it is not difficult to verify from the Codazzi equations that  the first normal spaces form a parallel subbundle of the normal bundle. This can also be derived from Theorem $2$ of \cite{dt2}. Then, it is a standard fact that the submanifold reduces codimension to two, i.e., it is a submanifold of codimension two of  a totally geodesic submanifold of the ambient space. Therefore, our main result provides a complete classification of all non-minimal submanifolds (of arbitrary codimension) of dimension at least four that attain equality in the inequality $(*)$ and whose first normal spaces have  dimension two everywhere. It is a very interesting problem to study the remaining cases.
   
   Notice that minimal submanifolds that attain equality in $(*)$ have necessarily first normal spaces of dimension two but may have arbitrary codimension. These submanifolds were considered in \cite{df2}. In particular, it was shown that complete examples must be Riemannian products $L^3\times \R^{n-3}$. Moreover, when the manifold is Kaehler a complete classification in terms of a Weierstrass-type representation was given. By composing such submanifolds with an inversion in Euclidean space one obtains non-minimal submanifolds that attain equality in $(*)$ whose first normal spaces have dimension three but do not reduce codimension.   

\section[The proofs]{The proofs}
  
 We first prove the converse of Theorem \ref{thm:main}. Let $f\colon\,M^n\to \R^{n+2}$ be an isometric immersion attaining equality everywhere in the inequality $(*)$ and free of minimal and umbilical points. We must prove that there exists  a  minimal surface $g\colon\,L^2\to \R^{n+2}$ and  a diffeomorphism $\psi\colon\,\Lambda_1\to M^n$, with $\Lambda_1$ defined as in the statement,  such that $\phi:=f\circ \psi$ is given by (\ref{eq:phi}).

By the result in \cite{ddvv}, there exist orthonormal tangent and normal frames $\{e_1,\ldots, e_n\}$ and $\{\eta, \zeta\}$, respectively, with respect to which the shape operators $A_\eta$ and $A_\zeta$ have the form (\ref{eq:shapeop}). Our assumption that $f$ is free of minimal and umbilical points is equivalent  to $\lambda$ and $\mu$ being nowhere vanishing, respectively. By Lemma 5.2 in \cite{ddvv}, we also have that $e_k(\lambda)=0$ and  $\nabla^\perp_{e_k}\eta=0$ for $k\geq 3$. Therefore $\lambda\eta$ is a {\em Dupin principal normal\/} of multiplicity $n-2$.   This means that the subspaces 
$$
E_\eta(x)=\{T\in T_xM\colon\,\alpha_f(T,X)=\lambda \<T,X\>\eta,\;\;
\mbox{for all}\,\,X\in T_xM\},
$$
where $\alpha_f$ is the second fundamental form of $f$ with values in the normal bundle, define a smooth distribution $E_\eta$ of rank $n-2$ satisfying
\be\label{eq:dupin}
T(\lambda)=0\,\,\,\mbox{and}\,\,\,\nabla_T^\perp \eta=0\,\,\,\,\mbox{for any}\,\,\,T\in E_\eta.
\ee
 In addition, since $\mu\neq 0$ we have that $\lambda\eta$ is {\em generic\/}, in the sense that $E_\eta=\ker(A_\eta-\lambda I).$
It is well-known that $E_\eta$ is an involutive distribution whose leaves are (mapped by $f$ into) round $(n-2)$-dimensional spheres in $\R^{n+2}$. 

    By the first equation in (\ref{eq:dupin}), the function $r=1/\lambda$ gives rise to a smooth function on the quotient space $L^2$ of leaves of $E_\eta$, which is also denoted by $r$. Let  $g\colon\,M^n\to \R^{n+2}$ be given by $g=f+r\eta$. From (\ref{eq:dupin}) we have  
$$g_*T=f_*T-rf_*A_\eta T=0,$$
hence $g$ also factors through a map on  $L^2$, still denoted by $g$.  By Proposition~$1$ in \cite{df} there exist a smooth unit vector field $\xi$ normal to $g$ and a diffeomorphism \mbox{$\psi\colon\,\Lambda_1\to M^n$}, where
$\Lambda_1$ is the unit bundle of the vector subbundle 
$\Lambda$ of the normal bundle of $g$ that is orthogonal to $\xi$,  such that 
$$
\eta(y,w)=(\eta\circ\psi)(y,w)=g_*\nabla r(y)+\rho(y)\xi(y)+\Omega(y) w
$$
for  $\rho,\Omega\in C^\infty(L)$ satisfying 
\be\label{eq:etanorm}
\|\nabla r\|^2+\rho^2+\Omega^2=1.
\ee
Moreover, we have
\be\label{eq:phi2}
\phi(y,w)=(f\circ \psi)(y,w)=g(y)-r(y)\eta(y,w).
\ee
We identify the tangent space $T_{(y,w)}\Lambda_1$ with the direct sum $T_yL\oplus \{w\}^\perp$, where $\{w\}^\perp$ denotes the orthogonal complement of $\spa\{w\}$ in $\Lambda(y)$, and write vectors $Y\in T_{(y,w)}\Lambda_1$ as $Y=(X,V)$ according to this decomposition. We also denote by  ${\cal V}$ the  corresponding {\em vertical\/} subbundle of $T\Lambda_1$, whose fiber ${\cal V}(y,w)$    at $(y,w)\in \Lambda_1$ is $\{w\}^\perp$.  Clearly, we have that $\psi_*{\cal V}=E_\eta$.
 
 Since $\lambda\eta$ is a  Dupin principal normal of $f$, the orthogonal complement $\eta^\perp$ of $\spa\{\eta\}$ in $T_f^\perp M$ is constant in $\R^{n+2}$ along $E_\eta$.  Therefore, if $\zeta$ is a unit vector field spanning $\eta^\perp$, then the map $\zeta\circ \psi$, which we also denote simply by $\zeta$, is constant along ${\cal V}$. Thus we may write 
\be\label{eq:zeta}
\zeta=g_*Z+a\xi\ee
for a smooth vector field $Z$ and $a\in C^\infty(L)$ satisfying
\be\label{eq:znorm}
\|Z\|^2+a^2=1.
\ee
Since $(\ker(A_\eta-\lambda I))^\perp$ has rank two everywhere, the function $a$ is nowhere vanishing. Otherwise $g_*Z$ would be somewhere normal to $\phi$, which would imply, by taking tangent
components for $X=Z$ in 
$$
\phi_*(X,0)=g_*X-\<\nabla r,X\>\eta-r\eta_*(X,0),
$$
 that $\psi_*(Z,0)\in \ker(A_\eta-\lambda I)$, a contradiction. \vspace{1ex}

 From now on we follow closely the proof of Theorem $1$ of \cite{dt}. 
\begin{lemma}\po\label{le:basic} 
The following facts hold:
\begin{itemize}
\item[$(i)$] $\rho=0$ and $L^2$ is orientable with a complex structure $J$ such that $JZ=\nabla r$,
\item[$(ii)$] $h=-r\zeta$ satisfies that $h_*=g_*\circ J$.
\end{itemize}
\end{lemma}

     Before proving Lemma \ref{le:basic}, let us see how it yields  the converse statement of the theorem. It follows from Lemma \ref{le:basic}-$(ii)$ that $g$ and $h$ are conjugate minimal surfaces. Moreover, from (\ref{eq:zeta}) we obtain  $-ra\xi=h^N$ and $-rg_*\nabla r=g_*Jh^T$, where $J$ is the complex structure on $L^2$ given by Lemma~\ref{le:basic}-$(i)$.  Since $\rho=0$ and $\|Z\|=\|\nabla r\|$ by Lemma \ref{le:basic}-$(i)$, it follows that $a^2=\Omega^2$ by (\ref{eq:etanorm}) and (\ref{eq:znorm}). Thus $|r\Omega |=\|h^N\|$, and then (\ref{eq:phi2}) reduces to  (\ref{eq:phi}). \vspace{1ex}

The proof of Lemma \ref{le:basic} will be given in several steps. We start with the following preliminary facts.

\begin{sublemma} {\hspace*{-1ex}}\textbf{. } \label{sub1}
 We have 
\be\label{bw20}
\<B_w Z,X\>=a\<\nabla_X^\perp w,\xi\>\,\,\,\,\,\,\mbox{for all $X\in TM$}\,\,\mbox{and}\,\,\,w\in \Lambda
\ee
and
\be\label{eq:hess}
\hess r(Z)-\frac{1}{r}Z
+B_\xi(a\nabla r-\rho Z)+a\nabla \rho=0,
\ee
where $B_w$ and $B_\xi$ denote the shape operators of $g$ in the normal directions $w$ and $\xi$, respectively.
\end{sublemma}
\proof
 From the fact that ${\cal V}=\ker A_\zeta$ we obtain
$$
0=\<\zeta_*(y,w)(X,0),\phi_*(y,w)(0,V)\>\,\,\,\mbox{for any}\,\,y\in L^2,\,w\in \Lambda_1(y)
\,\,\,\mbox{and}\,\,\,V\in w^\perp.
$$ 
Then (\ref{bw20}) follows by differentiating (\ref{eq:zeta}) and using that  
$\phi_*(0,V)=-r\Omega V$.
 Thus
\be\label{eq:zeta*}
\zeta_*(X,0)=g_*DX+\<K,X\>\xi,
\ee
where
\be\label{eq:dk}
DX=\nabla_X Z-aB_\xi X\;\;\;\;\mbox{and}\;\;\;\;\;K=\nabla a+B_\xi Z.
\ee
The orthogonality between $\eta$ and $\zeta$ yields
\be\label{eq:zr}
\<Z,\nabla r\>+a\rho=0.
\ee
Hence,
\begin{eqnarray}\label{eq:nxz}
\<\nabla_XZ,\nabla r\>\!\!&=&\!\!XZ(r)
-\<Z,\hess r(X)\>\nonumber\\
\!\!&=&\!\!-X(a)\rho -aX(\rho)-\<Z,\hess r(X)\>\nonumber\\
\!\!&=&\!\!-\<\rho\nabla a+a\nabla \rho+\hess r(Z),X\>.
\end{eqnarray}
It follows  from (\ref{eq:zeta*}), (\ref{eq:dk}) and (\ref{eq:nxz}) that
\be\label{eq:z1} 
\<\zeta_*(X,0),\eta\>=\<DX,\nabla r\>+\rho\<K,X\>=-\<\hess r(Z)+B_\xi(a\nabla r-\rho Z)+a\nabla \rho,X\>.
\ee
On the other hand,
\be\label{eq:z2}
\<\zeta,\phi_*(X,0)\>=X\<\zeta,g\>-\<\zeta_*(X,0),g-r\eta\>
=\<Z,X\>
+r\<\zeta_*(X,0),\eta\>,
\ee
thus (\ref{eq:hess})  follows from (\ref{eq:z1}), (\ref{eq:z2}) and the fact that $\zeta$ is normal to $\phi$.
\vspace{1,5ex}\qed

We now express in terms of $g$ the condition that the shape operators of $f$ are given by (\ref{eq:shapeop}). It is convenient to use the orthonormal frame $Y_1, Y_2, Y_j=e_j$ for $3\leq j\leq n$, with
$$
Y_1=\frac{1}{\sqrt{2}}(e_1+e_2),\;\;\;\;\;\; Y_2=\frac{1}{\sqrt{2}}(e_1-e_2)
$$ 
With respect to this frame, the matrices in  (\ref{eq:shapeop}) become
$$
A_\eta=\left[ \begin{array}{ccc}
\lambda +\mu & \;0 &\;0 \;\;\cdots\;\; 0\\ 
0 & \lambda -\mu&\;0 \;\;\cdots\;\; 0\\ 
0 & 0 &\;\lambda \;\;\cdots\;\; 0\\ 
\vdots & \vdots & \vdots \;\;\cdots\;\; \vdots\\ 
0 & 0 &\;0 \;\;\cdots\;\; \lambda\\
\end{array}
\right]\\;\;\;\;\;\;\;\;A_\zeta=\left[ 
\begin{array}{ccc}
0 & \;\mu &\;0 \;\;\cdots\;\; 0\\ 
\mu & \;0 &\;0 \;\;\cdots\;\; 0\\ 
0 & 0 &\;0 \;\;\cdots\;\; 0\\ 
\vdots & \vdots & \vdots \;\;\cdots\;\; \vdots\\ 
0 & 0 &\;0 \;\;\cdots\;\; 0\\
\end{array}
\right].
$$
Therefore,
\be\label{sff}  
\left\{ \begin{array}{l}
\eta_*Y_1=-(\lambda +\mu )\phi_*Y_1-\omega(Y_1)\zeta\vspace{1.5ex}\\
\eta_*Y_2=-(\lambda -\mu )\phi_*Y_2-\omega(Y_2)\zeta\vspace{1.5ex}\\
\zeta_*Y_i=-\mu \phi_*Y_j+\omega(Y_i)\eta,\,\,\,
1\leq i\neq j\leq 2,
\vspace{1.5ex}\\
\end{array} \right. 
\ee
where $\omega(Y)=\<\nabla^\perp_{Y}\zeta, \eta\>$ and
\be\label{eq:perp}
\<\phi_*Y_i(y,w),\phi_*(0,V)\>=0,\,\,1\leq i\leq 2,\,\,\,\mbox{for any}\,\,\,(y,w)\in \Lambda_1\,\,\,\mbox{and}\,\,\,V\in \{w\}^\perp.
\ee 
Write $Y_1=(X_1,V_1)$ and $Y_2=(X_2,V_2)$ according to the splitting $T_{(y,w)}\Lambda_1=T_yL\oplus \{w\}^\perp$. Using that 
$$
\eta_*(0,V)=\Omega
V,\,\,\,\,\phi_*(0,V)=-r\Omega
V\,\,\,\,\mbox{and}\,\,\,\zeta_*(0,V)=0,
$$
we obtain from (\ref{eq:perp}) that 
\be\label{eq:wperp}
(\eta_*(X_i,0))_{w^\perp}=-\Omega V_i,\,\,\,\,1\leq i\leq 2.
\ee
 Taking components in $H=g_*TL\oplus \spa\{\xi\}\oplus \spa\{w\}$ in (\ref{sff}) gives
\be\label{sff2}  
\left\{ \begin{array}{l}
r^2\mu(\eta_*(X_1,0))_{H}=\theta_1(g_*X_1-r_1\eta)+r\omega((X_1,0))\zeta\vspace{1.5ex}\\
r^2\mu(\eta_*(X_2,0))_{H}=\theta_2(-g_*X_2+r_2\eta)-r\omega((X_2,0))\zeta\vspace{1.5ex}\\
\zeta_*(X_i,0)=-\mu g_*X_j+\mu (r(\eta_*(X_j,0))_{H}+ r_j\eta)+\omega((X_i,0))\eta, \,\,\,1\leq i\neq j\leq 2,
 \end{array} \right. \ee
where 
$\theta_1=1+r\mu$,  $\theta_2=1-r\mu$ and $r_i=\<\nabla r, X_i\>$ for $1\leq i\leq 2$. 

We have 
$$ 
\eta_*(X,0)=g_*Q_wX
+\<T_w,X\>\xi+\<P_w,X\>w+(\eta_*(X,0))_{w^\perp},
$$
where
$$
\left\{ \begin{array}{l}
{\displaystyle Q_w=\hess r-\rho B_\xi -\Omega B_w}\label{eq:qw}\vspace{1ex}\\
{\displaystyle T_w=\nabla\rho +B_\xi\nabla
r+\frac{\Omega}{a}B_w Z} \label{eq:omega}\vspace{1ex}\\
{\displaystyle P_w=\nabla\Omega +B_w\nabla r
-\frac{\rho}{a}B_w Z}. \label{eq:alpha}\end{array} \right.
$$
It follows immediately from (\ref{eq:hess}) and (\ref{eq:z1}) 
that 
$$
\omega(X,0)=-\frac{1}{r}\<Z,X\>,
$$
in terms of the metric induced by $g$. 
Using this, the $w$-component of (\ref{sff2}) gives
\be\label{eq:wcomp}  
\left\{ \begin{array}{l}
r^2\mu\<P_w,X_1\>=-\theta_1\Omega r_1\vspace{1.5ex}\\
r^2\mu\<P_w,X_2\>=\theta_2\Omega r_2\vspace{1.5ex}\\
r^2\mu\<P_w,X_i\>=-r\mu\Omega r_i+\Omega\<Z,X_j\>,\,\,\,1\leq i\neq j\leq 2.
\end{array} \right. 
\ee
Replacing the first two equations into the last two yields
\be\label{eq:gradrZ}
r_1=-\<Z,X_2\> \;\;\;\mbox{and}\;\;\;\;\;r_2=\<Z,X_1\>.
\ee
Taking the tangent component to $g$ of (\ref{sff2}) and using (\ref{eq:gradrZ}) we obtain
\be\label{sff3}  \left\{ \begin{array}{l}
r^2\mu Q_wX_1-\theta_1SX_1+r_2Z=0\vspace{1.5ex}\\
r^2\mu Q_wX_2+\theta_2SX_2+r_1Z=0\vspace{1.5ex}\\
rDX_1+r\mu SX_2-r^2\mu Q_wX_2+r_2\nabla r=0\vspace{1.5ex}\\
rDX_2+r\mu SX_1-r^2\mu Q_wX_1-r_1\nabla r=0,\vspace{1.5ex}\\
 \end{array} \right. \ee
where we denoted
\be\label{eq:S}
 S=I-\<\nabla r,*\>\nabla r.
\ee

Finally, computing the $\xi$-component of (\ref{sff2}) gives
\be\label{eq:xicomp}  
\left\{ \begin{array}{l}
r^2\mu\<T_w,X_1\>=-\theta_1\rho r_1 -a\<Z,X_1\>\vspace{1.5ex}\\
r^2\mu\<T_w,X_2\>=\theta_2\rho r_2 +a\<Z,X_2\>\vspace{1.5ex}\\
r\<K, X_i\>=r^2\mu\<T_w,X_j\>+ r\mu\rho r_j
-\rho\<Z,X_i\>,\,\,\,1\leq i\neq j\leq 2.
\end{array} \right. 
\ee
We now use that
\be\label{eq:orth}
\delta_{ij}=\<\phi_*Y_i,\phi_*Y_j\>,\,\,\,1\leq i,j\leq 2.
\ee 
 From (\ref{eq:wperp}) we get
$$\<\phi_*(0,V_1),\phi_*(X_2,0)\>=\<-r\Omega V_1,-r\eta_*(X_2,0)\>=-r^2\Omega^2\<V_1,V_2\>=
\<\phi_*(0,V_2),\phi_*(X_1,0)\>.$$ On the other hand,
$$
\<\phi_*(0,V_1),\phi_*(0,V_2)\>=r^2\Omega^2\<V_1,V_2\>=
\<(\phi_*(X_1,0))_{w^\perp}, (\phi_*(X_2,0))_{w^\perp}\>.
$$
Hence (\ref{eq:orth}) reduces to $\delta_{ij}=\<(\phi_*(X_1,0))_{H},
(\phi_*(X_2,0))_{H}\>$, which gives
\be\label{eq:orth2}
\begin{array}{l}\delta_{ij}=\<X_i,X_j\>-r_ir_j-2r\<Q_wX_i,X_j\>+r^2(\<Q_wX_i,Q_wX_j\>\vspace{1ex}\\
\hspace*{5ex}+\<T_w,X_i\>\<T_w,X_j\>
+\<P_w,X_i\>\<P_w,X_j\>).\end{array}
\ee
Then, we argue exactly as in the proof of Sublemma $7$ in \cite{dt} to prove that
\be\label{eq:xinorm}\|X_1\|^2=r^2\mu^2+r_1^2+r_2^2=\|X_2\|^2.
\ee
 First, taking inner products of the first and second  equations in (\ref{sff3}) by $X_2$ and $-X_1$, respectively,
and adding them up, bearing in mind (\ref{eq:gradrZ}), yields
$$
\<X_1,X_2\>=0.
$$
On the other hand, we compute from the first two equations in (\ref{sff3}) that
\be\label{eq:qwuiuj}
\left\{\begin{array}{l}
r^2\mu\<Q_wX_1,X_1\>=\theta_1(\|X_1\|^2-r_1^2)-r_2^2
\vspace{1.5ex}\\
r^2\mu\<Q_wX_2,X_2\>=-\theta_2(\|X_2\|^2-r_2^2)+r_1^2
\vspace{1.5ex}\\
r\<Q_wX_1,X_2\>=-r_1r_2.\\
\end{array} \right.
\ee
Using (\ref{eq:etanorm}), (\ref{eq:znorm}), (\ref{eq:zr}) and the first two equations in (\ref{sff3}) we have
\be\label{eq:qwuiqwuj}
\hspace*{-2ex}\left\{\begin{array}{l}
\!\!r^4\mu^2\|Q_wX_1\|^2=\theta_1^2(\|X_1\|^2
-(1+\rho^2+\Omega^2)r_1^2)-2\theta_1(r_2^2+a\rho r_1r_2)
+(1-a^2)r_2^2\vspace{1.5ex}\\
\!\!r^4\mu^2\|Q_wX_2\|^2=\theta_2^2(\|X_2\|^2
-(1+\rho^2+\Omega^2)r_2^2)-2\theta_2(r_1^2-a\rho r_1r_2)
+(1-a^2)r_1^2\vspace{1.5ex}\\
\!\!r^4\mu^2\<Q_wX_1,Q_wX_2\>=(\theta_1\theta_2
(1+\rho^2+\Omega^2)-\theta_1-\theta_2-a^2+1)r_1r_2
-a\rho(\theta_1r_1^2-\theta_2r_2^2).\\
\end{array} \right.
\ee
 From (\ref{eq:gradrZ}) and (\ref{eq:xicomp}) we obtain
$$ r^2\mu\<T_w,X_1\>=-\theta_1\rho
r_1-ar_2\;\;\;\mbox{and}\;\;\; r^2\mu\<T_w,X_2\>=\theta_2\rho r_2-ar_1. $$
Thus,
\be\label{eq:twui2}  
\left\{ \begin{array}{l}
{\displaystyle r^4\mu^2\<T_w,X_1\>^2=\theta_1^2\rho^2r_1^2+a^2r_2^2
+2\theta_1a\rho r_1r_2}\vspace{1.5ex}\\{\displaystyle r^4\mu^2\<T_w,X_1\>\<T_w,X_2\>=
(a^2-\theta_1\theta_2\rho^2)r_1r_2
+\theta_1a\rho r_1^2}-\theta_2a\rho r_2^2
\vspace{1.5ex}\\
{\displaystyle r^4\mu^2\<T_w,X_2\>^2=\theta_2^2\rho^2r_2^2+a^2r_1^2-2\theta_2a\rho r_1r_2}.\end{array} \right. 
\ee
 From the first two equations in (\ref{eq:wcomp}) we get
\be\label{eq:pwui2}
\left\{ \begin{array}{l}
{\displaystyle r^4\mu^2\<P_w,X_1\>^2=\theta_1^2\Omega^2r_1^2}
\vspace{1.5ex}\\
{\displaystyle r^4\mu^2\<P_w,X_2\>^2=\theta_2^2\Omega^2r_2^2}
\vspace{1.5ex}\\
{\displaystyle r^4\mu^2\<P_w,X_1\>\<P_w,X_2\>=-\theta_1\theta_2
\Omega^2r_1r_2.}\end{array} \right.
\ee
Replacing (\ref{eq:qwuiuj}), (\ref{eq:qwuiqwuj}), (\ref{eq:twui2}) and (\ref{eq:pwui2}) into (\ref{eq:orth2})  we end up with (\ref{eq:xinorm}).

 It follows from (\ref{eq:gradrZ}) and (\ref{eq:xinorm}) that $\nabla r$ and $Z$ are orthogonal vector fields on $L^2$ with the same norm, thus there exists a complex structure $J$ on $L^2$ such that $JZ=\nabla r$. We conclude from (\ref{eq:zr}) that $\rho=0$, and the proof of $(i)$ is completed.\vspace{1ex}\\
We now prove $(ii)$.  Replacing the first two equations of (\ref{sff3}) into the last two gives
$$
\left\{ \begin{array}{l}
rDX_1+SX_2+r_1Z+r_2\nabla r=0\vspace{1.5ex}\\
rDX_2-SX_1+r_2Z-r_1\nabla r=0,\\
\end{array} \right.
$$
which can be written as
\be\label{eq:D}
rDX=-JX-\<\nabla r,X\>Z.
\ee
On the other hand, replacing the first two equations of (\ref{eq:xicomp}) into the last
two yields
\be\label{eq:rKui}
\left\{ \begin{array}{l}
r\<K,X_1\>=a\<Z,X_2\>
\vspace{1.5ex}\\
r\<K,X_2\>=-a\<Z,X_1\>.\\
\end{array} \right.
\ee
Taking (\ref{eq:gradrZ}) into account, the preceding equations reduce to
\be\label{eq:nablaar}
rB_\xi Z+\nabla(ar)=0.
\ee
 From (\ref{bw20}) we have
$$
\tilde{\nabla}_X\xi=-g_*B_\xi X+\nabla_X^\perp \xi\vspace{1ex}\\
=-g_*B_\xi X-\frac{1}{a}(\alpha_g(Z,X)-\<B_\xi Z,X\>\xi),
$$
where $\alpha_g$ denotes the second fundamental form of $g$. Hence,
$$
-ar\tilde{\nabla}_X\xi+r\<B_\xi Z,X\>\xi=arg_*B_\xi X+r\alpha_g(Z,X).
$$
In view of (\ref{eq:nablaar}) the left-hand-side is $\tilde{\nabla}_X(-ar\xi)$.
For the right-hand-side we have
\begin{eqnarray*}
arg_*B_\xi X+r\alpha_g(Z,X)\!\!&=&\!\!arg_*B_\xi X+r(\tilde{\nabla}_Xg_*Z-g_*\nabla_XZ)\vspace{1ex}\\
\!\!&=&\!\! g_*(arB_\xi X-r\nabla_XZ-X(r)Z)+\tilde{\nabla}_X(rg_*Z).
\end{eqnarray*}
Therefore, we obtain using (\ref{eq:D}) that
$$
h_*X=g_*(arB_\xi X-r\nabla_XZ-X(r)Z)=g_*(-rDX-X(r)Z)=g_*JX.
$$

   We now prove the direct statement of Theorem \ref{thm:main}. We need the following fact from~\cite{dt}.

\begin{proposition}\po\label{pr:conj}
Let $g\colon\,M^2\to \R^{n+2}$  be a simply connected oriented minimal surface with complex structure $J$
compatible with the orientation and let $h\colon\,M^2\to \R^{n+2}$ be a conjugate minimal surface such that $h_*=g_*\circ J$.
Then $r=\|h\|$ satisfies
$\|\nabla r\|\leq 1$ everywhere. Moreover, on the complement of the subset of isolated
points of $M^2$ where $a=\sqrt{1-\|\nabla r\|^2}$ vanishes, there
exists a  smooth unit normal vector field $\xi$ to $g$ 
such~that
$$
h=r(g_*J\nabla r-a \xi).
$$
Furthermore,  
\be\label{bw201} 
\<B_\delta J\nabla r,X\>+a\<\nabla_X^\perp\delta,\xi\>=0
\,\,\,\mbox{ for all}\,\,\,\delta\in \spa\{\xi\}^\perp
\ee 
and
\be\label{eq:bxi}
B_\xi=\frac{1}{ar}(r\mbox{Hess\,} r-S)\circ J,
\ee
where $S$ is given by (\ref{eq:S}).
\end{proposition}

Setting  $\eta(y,w)=g_*\nabla r(y) - a w$, we have by Proposition  \ref{pr:conj} that
$\phi=g-r\eta.$  
 From 
$$
\phi_*(X,0)=g_*X-\<\nabla r, X\>\eta-r\eta_*(X,0)\,\,\,\,\,\mbox{and}\,\,\,\phi_*(0,V)=-arV
$$
it follows that $\eta$ is a unit normal vector field to $\phi$. Moreover, since $\phi+r\eta=g$ does not depend on $w$, we have that ${A_\eta}|_{\cal V}=r^{-1}I$.
 
Let  $\zeta$  be defined by (\ref{eq:zeta}) with $Z=-J\nabla r$. 
Then $\zeta$  has unit length and is orthogonal to $\eta$. We obtain from (\ref{bw201}) that (\ref{eq:zeta*}) holds, hence we have (\ref{eq:z1})  with $\rho=0$, and  also  (\ref{eq:z2}).
 From (\ref{eq:bxi}) we get
$$
r\hess r(Z) -Z+arB_\xi\nabla r=0,
$$
which implies, using (\ref{eq:z1}) (with $\rho=0$)  and (\ref{eq:z2}), that $\zeta$ is normal to $\phi$. 
It also follows from (\ref{bw201}) that  
${A_\zeta}|_{\cal V}=0$.

 Therefore, to complete the proof it suffices to show that there exists an orthonormal frame $\{Y_1,Y_2\}$  on the open subset of regular points of $\phi$  (with respect to the metric induced by $\phi$) satisfying (\ref{sff})  and (\ref{eq:perp}).

For each $w\in \Lambda_1$, since $B_w$ and $B_\xi$ are traceless symmetric $2\times 2$ matrices, we have 
\be\label{eq:alpha2}(B_w+B_\xi J)^2=\alpha^2 I\ee  for some $\alpha\in \R$. We need only  prove the existence of the orthonormal frame $\{Y_1,Y_2\}$ on the complement of the subset with empty interior of points of $\Lambda_1$ where $\alpha$ is nonzero. At such a point, set  $\mu=-a/r^2\alpha$. Since 
$B_{w}+B_\xi J=\alpha R_w$ for some reflection 
$R_w$ by (\ref{eq:alpha2}), it follows using (\ref{eq:bxi}) that \be\label{eq:bw2}
B_{w}=\frac{1}{a}(\hess r-\frac{1}{r}S)-\frac{a}{r^2\mu}R_w.
\ee

     For each $w\in \Lambda_1$, let $\{\bar{X}_1, \bar{X}_2\}$ be the orthonormal basis of $TL$ (with respect to the metric induced by $g$) formed by eigenvectors of $R_w$, with $\bar{X}_1$ corresponding to the eigenvalue $+1$ and $\bar{X}_2=J\bar{X}_1$. Define 
$$X_i=\frac{r\mu}{a}\bar{X}_i,\,\,\, 1\leq i\leq 2,$$
and set 
$$
V_i=-\frac{1}{a}(\eta_*(X_i,0))_{w^\perp},\,\,\,\,1\leq i\leq 2.
$$
We claim that $\{Y_1,Y_2\}$ given by $Y_i=(X_i,V_i)$ is the desired 
orthonormal frame.

   It follows from $(\eta_*(X_i,0))_{w^\perp}=-aV_i$ that (\ref{eq:perp}) is satisfied. In particular,  in order to check that $\{Y_1,Y_2\}$ is an orthonormal frame  it suffices to verify  (\ref{eq:orth2}). 
It also follows from $(\eta_*(X_i,0))_{w^\perp}=-aV_i$ that  the $w^\perp$-components of both sides of all equations in (\ref{sff}) coincide.  Therefore, it suffices to prove that (\ref{sff2}), or equivalently, (\ref{eq:wcomp}), (\ref{sff3}) and (\ref{eq:xicomp}), holds for $X_1$ and $X_2$. 

  Since we have (\ref{eq:gradrZ}), because  $JX_1=X_2$ and $JZ=\nabla r$, system (\ref{eq:wcomp})
reduces to its first two equations. These are in turn equivalent to 
$$
rB_{w}\nabla r+\frac{a}{r\mu}R_w\nabla r+\nabla(ar)=0,
$$
which follows from (\ref{eq:bw2}).

Now, (\ref{eq:bw2}) also implies that
$$
rQ_w=S+\frac{a^2}{r\mu}R_w.
$$
Moreover,  from (\ref{eq:bxi}) we get (\ref{eq:D}), hence (\ref{sff3}) is satisfied.

 From (\ref{eq:bxi}) we obtain (\ref{eq:nablaar}), and hence (\ref{eq:rKui}). Moreover, (\ref{eq:bxi}) and (\ref{eq:bw2}) imply  that
$$
B_\xi\nabla r+B_{w} Z+\frac{a}{r^2\mu}R_w Z=0,
$$
thus 
(\ref{eq:xicomp}) is satisfied.

Finally, we now have (\ref{eq:qwuiuj}), (\ref{eq:qwuiqwuj}), (\ref{eq:twui2}) and (\ref{eq:pwui2}), hence  (\ref{eq:orth2}) follows by using that $\<X_1,X_2\>=0$ and $\|X_i\|=r\mu/{a}$ for $1\leq i\leq 2$.\vspace{1ex}\qed

%  It remains to verify that $\phi$ is singular at a point $(x_0,w_0)\in \Lambda_1$ where $\alpha$ vanishes. Since totally geodesic points of $g$ are isolated, there always exists a sequence $(x_n,w_n)$  in $\Lambda_1$ converging to $(x_0,w_0)$ so that $\alpha(x_n,w_n)\neq 0$. Let $\{\bar{X}_1, \bar{X}_2\}$ be the orthonormal basis of $TL$ (with respect to the metric induced by $g$)  previously defined, and let $X_i$, $1\leq i\leq 2$, be given as in (\ref{eq:xi}). Then  $\bar{X}_i=r\alpha X_i$, $1\leq i\leq 2$. Since $Y_i=(X_i,V_i)$ is a unit vector with respect to metric induced by $\phi$ and $\phi_*(0,V_i)=$ hence $(\bar{X}_i,0)$ has formed by eigenvectors of $R_w$, with $\bar{X}_1$ corresponding to the eigenvalue $+1$ and $\bar{X}_2=J\bar{X}_1$. Define 
%$$X_i=\frac{r\mu}{a}\bar{X}_i,\,\,\, 1\leq i\leq 2,$$Let $\bar{X}_i$ be 

\begin{remark}\po\label{reg} 
{\em It was shown in  \cite{df} that an isometric immersion $f\colon\,M^n\to \R^{n+p}$, $n\geq 4$, that carries a generic Dupin principal normal $\eta$ of multiplicity $n-2$  is a rotational submanifold over a surface whenever  $\mbox{trace}\,A_\eta\neq n \|\eta\|$ and $\mbox{trace}\,A_\eta$ is constant along the leaves of the corresponding eigendistribution. Our result shows that the assumption that $\mbox{trace}\,A_\eta\neq n \|\eta\|$ can not be removed.}
\end{remark}

For the proof of Theorem \ref{prop:minimas}, we first recall that, given an isometric immersion $f\colon\,M^n\to \R^N$ and an inversion with respect to a sphere of radius $d$ centered at $P_0\in \R^N$, the map
\be \label{eq:vbi}
{\cal P}\xi=\xi-2\frac{\<f-P_0,\xi\>}{\<f-P_0,f-P_0\>}(f-P_0)
\ee
is a vector bundle isometry between the normal bundles $T_f^\perp M$ and
$T_{{\cal I}\circ f}^\perp M$. Moreover,
the shape operators $A_\xi$ and $\tilde{A}_{{\cal P}\xi}$ of $f$ and ${\cal I}\circ f$ with respect to $\xi$ and ${\cal P}\xi$, respectively, are related by
\be\label{eq:sop}
\tilde{A}_{{\cal P}\xi}
=\frac{1}{d^2}\left(\<f-P_0,f-P_0\>A_\xi+2\<f-P_0,\xi\>I\right).
\ee
Similar results hold for an ``inversion"
$$
{\cal I}(P)=P_0-\frac{d^2}{\<P-P_0,P-P_0\>}(P-P_0),\;\;\;
P\neq P_0,
$$
in Lorentzian space $\Les^N$ with respect to a hyperbolic space 
$$
\Hy^{N-1}(P_0;d):=\{P\in \Les^N\,:\,\<P-P_0,P-P_0\>=-d^2\}
$$
 (see Lemma $15$ of \cite{dt}), with $d^2$ replaced by $-d^2$ in formula (\ref{eq:sop}). We also observe that a stereographic projection of $\Sf_d^{N}$ onto $\R^N$ can be regarded as the restriction to $\Sf_d^{N}$ of an inversion in $\R^{N+1}$ with respect to the sphere of radius $2d$ centered at $2de_{N+1}$. Similarly, a stereographic projection of $\Hy_d^N$ onto $B(0;2d)\subset \R^{N}$ can be viewed as the restriction to $\Hy_d^N$ of an inversion in $\Les^{N+1}$ with respect to $\Hy^{N}(-2de_{N+1};2d)$. In both cases, we regard $\R^N$ as the hyperplane through the origin orthogonal to $e_{N+1}$ in either $\R^{N+1}$ or $\Les^{N+1}$, respectively.

Let us denote by $\Q^{n+2}_{\epsilon}$ either $\R^{n+2}$, $\Sf_d^{n+2}$ or $\Hy_d^{n+2}$, according as $\epsilon =0$, $\epsilon =1$ or $\epsilon =-1$, respectively. Given an austere isometric immersion $f\colon\,M^n\to \Q^{n+2}_{\epsilon}$ attaining equality in the inequality $(*)$,  let
$\hat{J}$ be the complex structure on $T_f^\perp M$ determined by the
opposite orientation to that induced by the vector bundle isometry ${\cal
P}\colon\,T_f^\perp M\to  T_{{\cal I}\circ f}^\perp M$ given by (\ref{eq:vbi})  from the
orientation on $T_{{\cal I}\circ f}^\perp M$ defined by the orthonormal frame $\{\eta,\zeta\}$ as in (\ref{eq:shapeop}). Here  ${\cal I}$ is an inversion in $\R^{n+2}$ 
with respect to a 
hypersphere of radius $2d$ centered at the origin or a 
stereographic projection of 
$\Sf_d^{n+2}$ or $\Hy_d^{n+2}$ onto $\R^{n+2}$ or $B(0;2d)\subset \R^{n+2}$, according as $\epsilon =0$, $\epsilon =1$ or $\epsilon =-1$, respectively. Set also $\bar{\epsilon}=1$ if $\epsilon=1$ or $0$ and $\bar{\epsilon}=-1$ if ${\epsilon}=-1$.

\begin{proposition}\po\label{prop:pairs}  Let $f\colon\,M^n\to \Q^{n+2}_{\epsilon}$ be an austere isometric immersion  that attains equality in the inequality $(*)$. Then  the holomorphic representative $G$ of the minimal surface
    associated to   ${\cal I}\circ f$
%where ${\cal I}$ is an inversion in $\R^{n+2}$ 
%with respect to a 
%hypersphere centered at the origin or a 
%stereographic projection of 
%$\Sf_R^{n+2}$ or $\Hy_R^{n+2}$ onto $\R^{n+2}$ or $B(0;2R)\subset \R^{n+2}$, respectively, 
is given by
%$$
%G=\frac{R^2}{2\|f^N\|^2}({f^N}+i\hat{J}f^N)
%$$
%$$
%G=\epsilon 2Re_{n+3}+\frac{R^2}{2\|(f-\epsilon 
%2Re_{n+3})^N\|^2}((f-\epsilon 2Re_{n+3})^N+i\hat{J}(f-\epsilon 2Re_{n+3})^N),
%$$
\be\label{eq:if}
G=\epsilon 2de_{n+3}+2\bar{\epsilon}d^2\frac{(f-\epsilon 2de_{n+3})^N+i\hat{J}(f-\epsilon 2de_{n+3})^N}{\<(f-\epsilon 
2de_{n+3})^N,(f-\epsilon 
2de_{n+3})^N\>},
\ee
where $(f-\epsilon 2de_{n+3})^N$ denotes the normal component (in $\Q^{n+2}_{\epsilon}$) of the 
position vector  $f-\epsilon 2de_{n+3}$ in either $\R^{n+2}$, $\R^{n+3}$ or $\Les^{n+3}$, according as $\epsilon =0$, $\epsilon =1$ or $\epsilon =-1$, respectively. 
\end{proposition}
\proof Set $\bar{d}=2d$ and $P_0=\epsilon \bar{d}e_{n+3}$. Define
\be\label{eq:tildes}
\tilde{\zeta}=(\bar{\lambda}^2
+\bar{\nu}^2)^{-1/2}(\bar{\nu}{\cal P}\eta-
\bar{\lambda}{\cal P}\zeta)\;\;\;\;\;\mbox{and}\;\;\;\;\;
\tilde{\eta}=(\bar{\lambda}^2
+\bar{\nu}^2)^{-1/2}(\bar{\lambda}{\cal P}\eta+\bar{\nu}{\cal P}\zeta),
\ee
where ${\cal P}$ is given by (\ref{eq:vbi}), 
$$
\bar{\epsilon}\bar{d}^2\bar{\lambda}= 2\<f-P_0,\eta\>\,\,\,\,\,\,\mbox{and}\,\,\,\,\,\,
\bar{\epsilon}\bar{d}^2\bar{\nu}= 2\<f-P_0,\zeta\>.
$$
Using (\ref{eq:sop}),  we obtain that the shape operators $\tilde{A}_{\tilde{\eta}}$ and
$\tilde{A}_{\tilde{\zeta}}$ of
$\tilde{f}= {\cal I}\circ f$
are given as in (\ref{eq:shapeop}) with $\lambda$ and $\mu$ replaced, respectively, by
\be\label{eq:tildes2}
{\displaystyle \tilde{\lambda}=(\bar{\lambda}^2
+\bar{\nu}^2)^{1/2}\;\;\;\;\mbox{and}\;\;\;\;\tilde{\mu}=\frac{\<f-P_0,f-P_0\>}{\bar{\epsilon}\bar{d}^2}\mu}.
\ee
 The holomorphic curve $G={g}+i{h}$ associated to $\tilde{f}$ is given by 
\be\label{eq:til}
{g}=\tilde{f}+\tilde{r}\tilde{\eta}
=\tilde{f}+\frac{\bar{\lambda}{\cal P}\eta
+\bar{\nu}{\cal P}\zeta}{\bar{\lambda}^2+\bar{\nu}^2}
\;\;\;\;\mbox{and}\;\;\;\;
{h}=-\tilde{r}\tilde{\zeta}
=-\frac{\bar{\nu}{\cal P}\eta
-\bar{\lambda}{\cal P}\zeta}{\bar{\lambda}^2
+\bar{\nu}^2},
\ee
where $\tilde{r}=1/\tilde{\lambda}$. We have
\be\label{eq:r4}
\bar{d}^4(\bar{\lambda}^2+\bar{\nu}^2)=4(\<f-P_0,\zeta\>^2
+\<f-P_0,\eta\>^2)=4\<(f-P_0)^N,(f-P_0)^N\>.
\ee
On the other hand, from
$${\cal P}\eta=\eta-2\frac{\<f-P_0,\eta\>}{\<f-P_0,f-P_0\>}(f-P_0)
\;\;\;\;
\mbox{and}\;\;\;\;
{\cal P}\zeta=\zeta-2\frac{\<f-P_0,\zeta\>}{\<f-P_0,f-P_0\>}(f-P_0)
$$
we obtain
\be\label{eq:peta}
 \bar{\epsilon}\bar{d}^2(\bar{\lambda}{\cal P}\eta
+\bar{\nu}{\cal P}\zeta)
=2(f-P_0)^N-4\frac{\<(f-P_0)^N,(f-P_0)^N\>}{\<f-P_0,f-P_0\>}(f-P_0)
\ee
and
\be\label{eq:pzeta}
\bar{\epsilon}\bar{d}^2(\bar{\nu}{\cal P}\eta-\bar{\lambda}{\cal P}\zeta)
=2\<f-P_0,\zeta\>\eta-2\<f-P_0,\eta\>\zeta=-2\hat{J}(f-P_0)^N.
\ee
  Then (\ref{eq:if}) follows from (\ref{eq:til}), (\ref{eq:r4}),  (\ref{eq:peta}) and
(\ref{eq:pzeta}).\vspace{2ex}\qed

 \noindent {\em Proof of Theorem \ref{prop:minimas}:\/} It follows from Proposition \ref{prop:pairs} that the holomorphic curve $G={g}+i{h}$ associated to ${\cal I}\circ f$ satisfies 
\be\label{eq:gh}\<g,h\>=0\,\,\,\,\,\mbox{and}\,\,\,\,\,\<g-P_0,g-P_0\>=\<h,h\>,\ee
with $P_0=\epsilon 2de_{n+3}$. Hence $G$ takes values in ${\cal H}_{\epsilon 4 d^2}$.

     Conversely, assume that the holomorphic curve $G={g}+i{h}$ associated to $f$
satisfies (\ref{eq:gh}). We claim that ${\cal I}\circ f$ is an austere isometric immersion into $\Q^{n+2}_{\epsilon}$. 

 Define
$
\tilde{\zeta}$ and $\tilde{\eta}$ as in (\ref{eq:tildes}), 
where ${\cal P}$ is the vector bundle isometry between
$T_f^\perp M$ and $T_{{\cal I}\circ f}^\perp M$ given by (\ref{eq:vbi}) and $\bar{\nu}, \bar{\lambda}$ are now given by
$$
\bar{\epsilon}\bar{d}^2\bar{\nu}=
2\<f-P_0,\zeta\>\;\;\;\;\;\mbox{and}\;\;\;\;\;
\bar{\epsilon}\bar{d}^2\bar{\lambda}= \lambda\<f-P_0,f-P_0\>+ 2\<f-P_0,\eta\>.
$$
As before, we obtain using (\ref{eq:sop})  that the shape operators
$\tilde{A}_{\tilde{\eta}}$ and $\tilde{A}_{\tilde{\zeta}}$ of $
{\cal I}\circ f$ are given as in (\ref{eq:shapeop}) with $\lambda$
and $\mu$ replaced, respectively, by $\tilde{\lambda}$ and $\tilde{\mu}$ given by (\ref{eq:tildes2}).
Since  (\ref{eq:gh}) holds, using that
${\displaystyle h=-\frac{1}{\lambda}\zeta}$ and ${\displaystyle
g-P_0=f-P_0+\frac{1}{\lambda}\eta}$  we obtain
$$
\<f-P_0,\zeta\>=\<g-P_0-\frac{1}{\lambda}\eta,\zeta\>=0$$ and
\begin{eqnarray*}
\frac{-2}{\lambda}\<f-P_0,\eta\>&=&\<f-P_0,f-P_0\>+\frac{1}{\lambda^2}-\<g-P_0,g-P_0\>\\&=&\<f-P_0,f-P_0\>+
\frac{1}{\lambda^2}-\<h,h\>=\<f-P_0,f-P_0\>.
\end{eqnarray*} Thus,
$\bar{\nu}=0=\bar{\lambda}$, and hence $\tilde{\lambda}=0$. Therefore ${\cal I}\circ f$ is  austere.
 \vspace{2ex}\qed

  The proof of Theorem \ref{thm:pairs} is  similar to that of Theorem $2$ of \cite{dt} and will be omitted.

\vspace{.2in} {\renewcommand{\baselinestretch}{1}

\hspace*{-20ex}\begin{tabbing} \indent\= IMPA -- Estrada Dona Castorina, 110
\indent\indent\=  Universidade Federal de São Carlos\\
\> 22460-320 -- Rio de Janeiro -- Brazil  \>
13565-905 -- São Carlos -- Brazil \\
\> E-mail: marcos@impa.br \> E-mail: tojeiro@dm.ufscar.br
\end{tabbing}}

\begin{thebibliography}{lbllll}

\bibitem{br} R. Bryant,  {\em Some remarks on the geometry of austere manifolds.\/} Bol. Soc. Bras. Mat. {\bf 21} (1991), 122--157.

\bibitem{ch}  B.Y. Chen, {\em Mean curvature and shape operators of isometric immersions into real-space-forms.\/}  Glasgow Math. J. {\bf 38}  (1996), 87--97.

\bibitem{cl}  T. Choi and Z. Lu, {\em On the DDVV conjecture and the comass in calibrated geometry.\/} Preprint (available at arXiv:math.DG/0610709).

\bibitem{df} M. Dajczer and L. Florit, {\em On Chen's Basic Equality.\/} Illinois J. Math.   {\bf 42} (1) (1998), 97--106.

\bibitem{df2} M. Dajczer and L. Florit, {\em A class of austere submanifolds.\/} Illinois J. Math. {\bf 45}  (2001), 735--755.

\bibitem{dt} M. Dajczer and R. Tojeiro, {\em All superconformal surfaces in $\R^4$ in terms of minimal surfaces.\/} Preprint (available at arXiv:math.DG/0710.5317).

\bibitem{dt2} M. Dajczer and R. Tojeiro, {\em  Submanifolds with nonparallel first normal bundle.\/}
Can. Math. Bull. {\bf 37} (1994), 330--337.

\bibitem{dfv}  F. Dillen, J. Fastenakels, J. Van der Veken, {\em  Remarks on an inequality involving the normal scalar curvature.\/} Preprint (available at arXiv:math.DG/0610721).

\bibitem{gt}   J. Ge and Z. Tang, {\em A proof of the DDVV conjecture and its equality case.\/} Preprint (available at arXiv:math.DG/0801.0650).

\bibitem{gr} I. Guadalupe and L. Rodríguez, {\em Normal curvature of surfaces in space forms.\/}  Pacific J. Math.   {\bf 106}  (1983), 95--103.

\bibitem{ho} D. Hoffman and R. Osserman, {\it The geometry of the generalized Gauss map.}
Mem. Amer. Math. Soc.  {\bf 28}  (1980), no. 236.

\bibitem{la} B. Lawson Jr, Lectures on minimal submanifolds. Mathematics Lecture Series, vol. $9$, Publish or Perish Inc., Wilmington, Del., 1980.

\bibitem{cl3}   Z. Lu, {\em On the DDVV conjecture and the comass in calibrated geometry $(II)$.\/} Preprint (available at arXiv:math.DG/0708.2921).

\bibitem{cl2}   Z. Lu, {\em Proof of the normal scalar curvature conjecture.\/} Preprint (available at arXiv:math.DG/0711.3510).

\bibitem{ddvv} P.J. de Smet, F. Dillen, L. Verstraelen and L. Vrancken, {\em  A pointwise inequality in submanifold theory.\/}  Arch. Math. (Brno)   {\bf 35}  (1999), 115--128.

\end{thebibliography}
\end{document}